\documentclass{article}
\usepackage{amssymb,latexsym}
\usepackage{amsmath}
\usepackage{graphics}
\usepackage{graphicx}
\newtheorem{theorem}{Theorem}

\newtheorem{definition}{Definition}

\newenvironment{proof}[1][Proof]{\noindent\textbf{#1.} }{\ \rule{0.5em}{0.5em}}
\numberwithin{equation}{section} \numberwithin{theorem}{section}
\numberwithin{corollary}{section}
\numberwithin{definition}{section} \numberwithin{result}{section}

\begin{document}
\title{A Note an Admissible Mannheim Curves in Galilean Space $G_3$}
\author{S. ERSOY, M. AKY\.{I}\~{G}\.{I}T, M. TOSUN  \\
sersoy@sakarya.edu.tr, maky\.{i}g\.{i}t@sakarya.edu.tr, tosun@sakarya.edu.tr \\
Department of Mathematics, Faculty of Arts and Sciences,\\
Sakarya University, Sakarya, 54187 TURKEY}

\maketitle
\begin{abstract}
The aim of this paper is to study the Mannheim partner curves in
three dimensional Galilean space $G_3$. Some well known theorems
are obtained related to Mannheim curves.


\end{abstract}

\section{Introduction}\label{S:intro}
In 1850 J. Bertrand defined Bertrand curves, with the help of this
definition Liu and Wang called the Mannheim pair which the
principal normal vector of the first curve coincides with the
binormal vector of the second curve and they obtained the
necessary and sufficient conditions between the curvature and the
torsions for a curve to be the Mannheim partner curves,
\cite{Liu}. Moreover, Orbay and Kasap examined on Mannheim curves
in $\mathbb{E}^3$, \cite{Orb}. The geometry of the Galilean space
$G_3$ has been treated in detail in R\"{o}schel's habilitation,
\cite{Ros}. Futhermore, Kamenarovic and Sipus studied about
Galilean space, \cite{Sip, Kam}. The properties of the curves in
the Galilean space  are studied in \cite{Pav}, \cite{ogr1},
\cite{ogr2}. In this paper, we gave some theorems and relations
about the curvatures and torsions of admissible Mannheim curves in
3-dimensional Galilean space $G_3$.

\section{Preliminaries}\label{Sec2}
The Galilean space $G_3$ is a Cayley-Klein space equipped with the
projective metric of signature (0,0,+,+), as in \cite{Mol}. The
absolute figure of the Galilean Geometry consist of an ordered
triple $\{ w,f,I\}$, where $w$ is the ideal (absolute) plane, $f$
is the line (absolute line) in $w$ and $I$ is the fixed elliptic
involution of points of $f$, \cite{Sip}. In the non-homogeneous
coordinates the similarity group $H_8$ has the form
\begin{equation}\label{2.1}
\begin{array}{l}
 \overline x  = a_{11}  + a_{12} x \\
 \overline y  = a_{21}  + a_{22} x + a_{23} y\cos \varphi  + a_{23} z\sin \varphi  \\
 \overline z  = a_{31}  + a_{32} x - a_{23} y\sin \varphi  + a_{23} z\cos \varphi  \\
 \end{array}
\end{equation}
where $a_{ij}$ and $\varphi$ are real numbers, \cite{Pav}.\\
In what follows the coefficients $a_{12}$ and $a_{23}$ will play
the special role. In particular, for $a_{12}  = a_{23}  = 1$,
(\ref{2.1}) defines the group $B_6  \subset H_8$ of isometries of
Galilean space $G_3$.\\
In $G_3$ there are four classes of lines:\\
\textit{i)} (proper) non-isotropic lines- they don't meet the
absolute line $f$.\\
 \textit{ii)} (proper) isotropic lines- lines that don't belong to the plane $w$ but meet the absolute line $f$.\\
 \textit{iii)} unproper non-isotropic lines-all lines of $w$ but $f$.\\
 \textit{iv)} the absolute line $f$.\\
Planes $x =$constant are Euclidean and so is the plane $w$. Other
planes are isotropic, \cite{Kam}.\\
Galilean scalar product can be written as
\[ < {u_1},{u_2} >  = \left\{ \begin{array}{l}
 \,\,\,\,\,\,\,{x_1}{x_2}\,\;\;\;\;\;\;,\;\;\;if\;\;{x_1} \ne 0\; \vee \;{x_2} \ne 0 \\
 {y_1}{y_2} + {z_1}{z_2}\;,\;\;if\;\;{x_1} = 0\; \wedge \;{x_2} = 0 \\
 \end{array} \right.\]
where ${u_1} = ({x_1},{y_1},{z_1})$ and ${u_2} =
({x_2},{y_2},{z_2})$. It leaves invariant the Galilean norm of the
vector $u = (x,y,z)$ defined by
\[\left\| u \right\| = \left\{ \begin{array}{l}
 \,\,\,\,\,\,\,\,\,x\;\;\;\;\;\;\;\;\;,\;x \ne 0 \\
 \sqrt {{y^2} + {z^2}} \;\;,\;x = 0. \\
 \end{array} \right.\]
Let $\alpha $ be a curve given in the coordinate form
\begin{equation}\label{2.2}
\begin{array}{l}
 \alpha :I \to {G_3},\;\;\;\;\;\;I \subset \mathbb{R}  \\
 \;\;\;\;\;\;t \to \alpha (t) = (x(t),y(t),z(t)) \\
 \end{array}
\end{equation}
where $x(t),y(t),z(t) \in {C^3}$ and $t$ is a real interval. If
$x'(t) \ne 0$, then the curve $\alpha $ is called admissible
curve.\\
Let $\alpha $ be an admissible curve in $G_3$, is parameterized by
arc length $s$, is given by
\[\alpha (s) = (s,y(s),z(s))\]
where the curvature ${\kappa}(s)$ and the torsion ${\tau }(s)$ are
\begin{equation}\label{2.3}
\kappa (s) = \sqrt {y'{'^2}(s) + z'{'^2}(s)\;}
\;\,\,\,\,\,\,\,{\rm{and }}\,\,\,{\rm{ }}\,\,\tau (s) =
\frac{{\det \left[ {\alpha '\left( s \right),\alpha ''\left( s
\right),\alpha '''\left( s \right)} \right]}}{{{\kappa ^2}(s)}},
\end{equation}
respectively. The associated moving Frenet frame is
\begin{equation}\label{2.4}
\begin{array}{l}
 {T}(s) = \alpha '(s) = (1,y'(s),z'(s)) \\
 {N}(s) = \frac{1}{{{\kappa}(s)}}\alpha '(s) = \frac{1}{{{\kappa}(s)}}(0,y''(s),z''(s)) \\
 {B}(s) = \frac{1}{{{\kappa}(s)}}(0, - z''(s),y''(s)). \\
 \end{array}
\end{equation}
Here ${T},{N}$ and ${B}$ are called the tangent vector, principal
normal vector and binormal vector fields of the curve $\alpha$,
respectively. Then for the curve $\alpha$, the following Frenet
equations are given by
\begin{equation}\label{2.5}
\begin{array}{l}
 {T}'(s) = {\kappa}(s){N}(s) \\
 {N}'(s) = {\tau }(s){B}(s) \\
 {B}'(s) =  - {\tau}(s){N}(s) \\
 \end{array}
\end{equation}
where ${T},{N},{B}$ are mutually orthogonal vectors, \cite{Kam}.

\section{Admissible Mannheim Curves in Galilean Space $G_3$}\label{Sec3}
In this section, we defined the admissible Mannheim curve and gave
some theorems related to these curves in $G_3$.
\begin{definition}\label{T3.1.}
Let $\alpha$ and $\alpha ^*$ be an admissible curves with the
Frenet frames along $\{ {T},{N},{B}\}$ and $ {\kern 1pt} \{ T^*,N
^*,B^*\} $, respectively. The curvature and torsion of $\alpha$
and $\alpha ^*$, respectively, ${\kappa}(s),{\tau}(s)$ and
$\kappa^*(s),\tau*(s)$ never vanish for all $s \in I$ in $G_3$. If
the principal normal vector field ${N}$ of $\alpha$ coincidence
with the binormal vector field $B^*$ of $\alpha ^*$ at the
corresponding points of the admissible curves $\alpha$ and $\alpha
^*$. Then $\alpha$ is called an admissible Mannheim curve and
$\alpha ^*$ is an admissible Mannheim mate of $\alpha$. Thus, for
all $s \in I$
\begin{equation}\label{3.1}
{\alpha ^*}(s) = \alpha (s) + \lambda (s){N}(s).
\end{equation}
\end{definition}
The mate of an admissible Mannheim curve is denoted by $(\alpha
,{\alpha ^*})$, see Figure 1.\\
\begin{center}
\hfil\scalebox{1}{\includegraphics{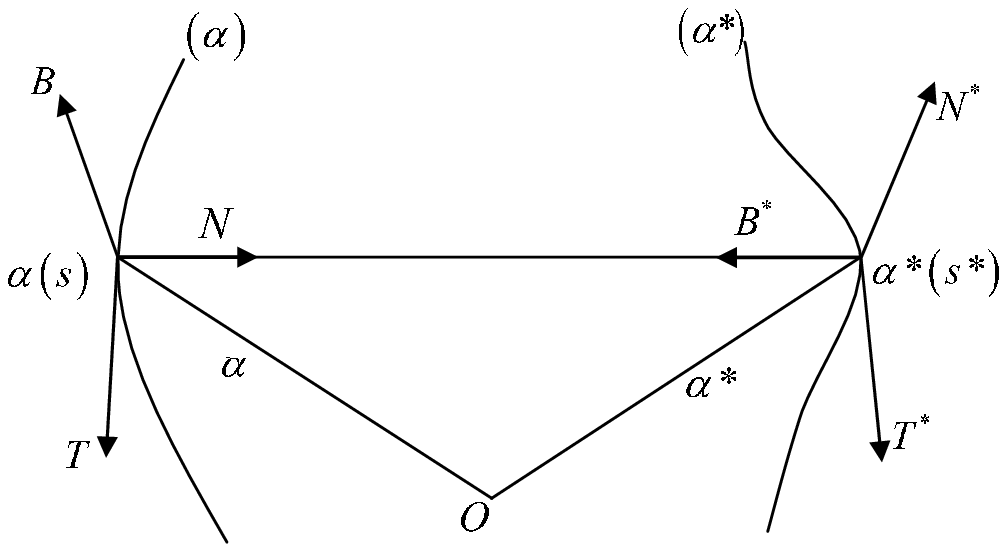}}\hfil\\
\scriptsize{Figure 1. The admissible Mannheim partner curves}\\
\end{center}

\begin{theorem}\label{T3.1.}
Let $(\alpha ,{\alpha ^*})$ be a mate of admissible Mannheim pair
in $G_3$. Then function $\lambda$ is constant on $I$ and defined
by equation (\ref{3.1}).
\end {theorem}

\begin{proof}
Let $\alpha$ be an admissible Mannheim curve in $G_3$ and
$\alpha^*$ be an admissible Mannheim mate of $\alpha$. Let the
pair of $\alpha (s)$ and ${\alpha ^*}(s)$ be of corresponding
points of $\alpha$ and $\alpha ^*$. Then the curve ${\alpha
^*}(s)$ is given by (\ref{3.1}). Differentiating (\ref{3.1}) with
respect to $s$ and using Frenet equations,
\begin{equation}\label{3.2}
T^*\frac{{d{s^*}}}{{ds}} = {T} + \lambda '{N} + \lambda {\tau}{B}
\end{equation}
is obtained.\\
Here and here after prime denotes the derivative with respect to
$s$. Since ${N}$ is coincident with $B^*$ in the same direction,
we have
\begin{equation}\label{3.3}
\lambda '(s) = 0,
\end{equation}
that is, $\lambda$ is constant. This theorem proves that the
distance between the curve $\alpha $ and its Mannheim mate
$\alpha^*$ is constant at the corresponding points of them. It is
notable that if $\alpha^*$ is an admissible Mannheim mate of
$\alpha$. Then $\alpha$ is also Mannheim mate of $\alpha^*$
because the relationship obtained in theorem between a curve and
its Mannheim mate is reciprocal one.
\end{proof}\\
\begin{theorem}\label{T3.2.}
Let $\alpha$ be an admissible curve with arc length parameter $s$.
$\alpha$ is an admissible Mannheim curve if and only if the
torsion ${\tau}$ of $\alpha$ is constant.
\end {theorem}
\begin{proof}
Let $(\alpha ,{\alpha ^*})$ be a mate of an admissible Mannheim
curves, then there exists the relation
\begin{equation}\label{3.4}
\begin{array}{l}
 {T}(s) = \cos \theta T^*(s) + \sin \theta N^*(s) \\
 {B}(s) =  - \sin \theta T^*(s) + \cos \theta N^*(s) \\
 \end{array}
\end{equation}
and
\begin{equation}\label{3.5}
\begin{array}{l}
 T^*(s) = \cos \theta {T}(s) - \sin \theta {B}(s) \\
 N^*(s) = \sin \theta {T}(s) + \cos \theta {B}(s) \\
 \end{array}
\end{equation}
where $\theta$ is the angle between ${T}$ and ${T}^*$ at the
corresponding points of $\alpha (s)$ and $\alpha^*
(s)$, (see Figure 1).\\
By differentiating (\ref{3.5}) with respect to $s$, we get
\begin{equation}\label{3.6}
\begin{array}{l}
\tau^*B^*\frac{{d{s^*}}}{{ds}} = \frac{{d(\sin \theta )}}{{ds}}{T}
+ \sin \theta {\kappa}{N} + \cos \theta {\tau}{N} + \frac{{d(\cos
\theta )}}{{ds}}{B}.
 \end{array}
\end{equation}
Since the principal normal vector field ${N}$ of the curve
$\alpha$ and the binormal vector field $B^*$ of its Mannheim mate
curve, then it can be seen that $\theta$ is a
constant angle.\\
If we equations (\ref{3.2}) and (\ref{3.5}) is considered, then
\begin{equation}\label{3.7}
\begin{array}{l}
\lambda {\tau}\cot \theta  = 1
 \end{array}
\end{equation}
is obtained. According to Theorem 3.1 and constant angle $\theta$,
$u = \lambda \cot \theta$ is constant. Then from equation
(\ref{3.7}), ${\tau} = \frac{1}{u}$ is constant, too.\\
Hence the proof is completed.
\end{proof}\\
\begin{theorem}\label{T3.3.}
(Schell's Theorem). Let $(\alpha ,{\alpha ^*})$ be a mate of an
admissible Mannheim curves with torsions ${\tau}$ and ${\tau^* }$,
respectively. The product of torsions ${\tau}$ and ${\tau^* }$ is
constant at the corresponding points $\alpha (s)$ and $\alpha^*
(s)$.
\end {theorem}
\begin{proof}
Since $\alpha$ is an admissible Mannheim mate of ${\alpha ^*}$,
then (\ref{3.1}) also can be given by
\begin{equation}\label{3.8}
\begin{array}{l}
\alpha  = {\alpha ^*} - \lambda B^*.
 \end{array}
\end{equation}
By taking differentiation of last equation and using equation
(\ref{3.4}),
\begin{equation}\label{3.9}
\begin{array}{l}
\tau^* = \frac{1}{\lambda }\tan \theta
 \end{array}
\end{equation}
can be given. By the helps of (\ref{3.7}), the equation below is
obtained easily;
\begin{equation}\label{3.10}
\begin{array}{l}
{\tau}\tau^* = \frac{{{{\tan }^2}\theta }}{{{\lambda ^2}}}
=constant
 \end{array}
\end{equation}
This completes the proof.
\end{proof}\\

\begin{theorem}\label{T3.4.}
Let $(\alpha ,{\alpha ^*})$ be an admissible Mannheim mate with
curvatures ${\kappa },\kappa^*$ and torsions ${\tau},\tau^*$ of
$\alpha$ and ${\alpha ^*}$, respectively. Then their curvatures
and torsions satisfy the following relations
\[\begin{array}{l}
 i)\;\;\kappa^* =  - \frac{{d\theta }}{{d{s^*}}} \\
 ii)\;\;{\kappa} = \tau ^*\frac{{d{s^*}}}{{ds}}\sin \theta  \\
 iii)\;{\tau} =  - \tau^*\frac{{d{s^*}}}{{ds}}\cos \theta . \\
 \end{array}\]
\end {theorem}

\begin{proof}
\textit{i)} Let us consideration equation (\ref{3.4}), then we
have
\begin{equation}\label{3.11}
\begin{array}{l}
 < {T},T^* >  = \cos \theta .
 \end{array}
\end{equation}
By differentiating last equation with respect to $s^*$ and using
the Frenet equations of $\alpha$ and ${\alpha ^*}$, we reach
\begin{equation}\label{3.12}
\begin{array}{l}
 < {\kappa}(s){N}(s)\frac{{ds}}{{d{s^*}}},T^*(s) >  +  < {T}(s),\kappa ^*(s)N^*(s) >  =  - \sin \theta \frac{{d\theta }}{{d{s^*}}}.
 \end{array}
\end{equation}
Since the principal normal ${N}$ of $\alpha$ and binormal $B^*$ of
${\alpha ^*}$ are linearly dependent. By considering equations
(\ref{3.4}) and (\ref{3.12}), we reach
\begin{equation}\label{3.13}
\begin{array}{l}
\kappa^*(s) =  - \frac{{d\theta }}{{d{s^*}}}.
 \end{array}
\end{equation}
If we take into consideration $< \,T,{B^*}\, >$, $< \,B,{B^*}
>$ scalar products and (\ref{2.5}), (\ref{3.4}), (\ref{3.5}) equations, then we can
easily prove $\textit{ii)}$ and $\textit{iii)}$ items of the
theorem, respectively.\\
The relations given in $\textit{ii)}$ and $\textit{iii)}$ of the
last theorem, we obtain
\[\frac{{{\kappa}}}{{{\tau}}} =  - \tan \theta=constant. \]
\end{proof}\\

\end {document}